\documentclass[11pt]{article}
\usepackage{amsmath}

\newcommand{\Tr}{\mathrm{Tr}\,}
\newcommand{\rhoe}{\rho^{e}}
\newcommand{\rhoh}{\rho^{h}}
\newcommand{\rhozero}{\rho^{(0)}}
\newcommand{\Xzero}{X^{(0)}}
\newcommand{\XC}{X^{(C)}}
\newcommand{\XiC}{\Xi^{(C)}}
\newcommand{\supz}{^{(0)}}
\newcommand{\supC}{^{(C)}}
\newcommand{\wt}{\widetilde}
\newcommand{\thetaz}{\theta^{(0)}}
\newcommand{\Thetaz}{\Theta^{(0)}}
\newcommand{\xiz}{\xi^{(0)}}
\newcommand{\Thetah}{\Theta^h}
\newcommand{\thetah}{\theta^h}
\newcommand{\wh}{w^h}
\newcommand{\wz}{w^{(0)}}
\newcommand{\diag}{\mathrm{diag} \,}
\newtheorem{definition}{Definition}[section]
\newtheorem{theorem}[definition]{Theorem}
\newtheorem{lemma}[definition]{Lemma}

\newtheorem{corollary}[definition]{Corollary}
\allowdisplaybreaks
\title{
All the trajectories of an extended averaged Hebbian learning equation
on the quantum state space are the e-geodesics
}
\author{
Yoshio Uwano
\\
\\
Department of Mathematics, 
Kyoto Pharmaceutical University, 
\\
Misasagi Nakauchi-cho 5, Yamashina-ku, Kyoto, 607-8414, Japan
\\
\\
e-mail:  uwano@mb.kyoto-phu.ac.jp
}
\date{}
\begin{document}
\maketitle
\begin{abstract}
%\noindent
In this paper, two families of trajectories on the quantum state space (QSS)
originating from a synaptic-neuron model and from quantum information
geometry meet together. The extended averaged Hebbian learning equation
(EAHLE) on the QSS developed by the author and Yuya \cite{UY1} from
a Hebbian synaptic-neuron model is studied from a quantum-information-geometric
point of view. It is shown that all the trajectories of  the EAHLE are the
e-geodesics, the autoparallel curves with respect to the exponential-type
parallel transport, on the QSS. As a secondary outcome, an explicit representation
of solution of the averaged Hebbian learning equation, the origin of the EAHLE,
is derived from that of the e-geodesics on the QSS.
\end{abstract}
%
%
%%%%%%%%%%%%%%%%%%%%%%%%%%%
\section{Introduction}
Quantum computing and quantum information have been well-known to be highly
developing fields in which a number of disciplines such as quantum mechanics,
mathematics, communication, information, statistics, control, optimization,
etc.\ are crossing over (see \cite{NC} for their history, for example).   
In this paper, a pair of interesting mathematical objects having different origins
meet together both of which are described on the quantum state space (QSS),
the space of regular density matrices of an arbitrarily fixed degree:
One is the extended averaged Hebbian learning equation (EAHLE) on the QSS
\cite{UY1} originating from one of the Hebbian synaptic-neuron learning
models and the other is the e-geodesics arising naturally from quantum
information geometry.
\par
In order to describe the motive of this paper, a brief history of the series of
the author's works \cite{UY1,UHI,UY2,U} is given below. As a departure point
of the series of the author's works, Grover's quantum search algorithm \cite{G}
is worth touched on, which is well-known to be one of the milestones on the
road of quantum computation \cite{NC}: For a large number, say $N=2^n$,
of randomly sorted data, the complexity of Grover's algorithm is of
$O(\sqrt{N})$, which is lower than the theoretical boundary, $O(N)$, of any
non-quantum searches. On Grover's algorithm, Miyake and Wadati \cite{MW}
made a pioneering geometric study saying that the search sequence is on a
geodesic on the $2^{n+1}-1$ dimensional unit sphere, $S^{2^{n+1}-1}$, of
$n$-qubit states and that the projection of the search sequence on the complex
projective space ${\bf C}P^{2^{n}-1}$ is also on a geodesic. Motivated by
Miyake and Wadati \cite{MW}, the author made a geometric study, with Hino
and Ishiwatari, on a Grover-type search for an ordered tuple of multi-qubits
\cite{UHI}: As a rigorous analogue to the projection applied to Grover's search
sequence, a projection map from the space of ordered tuples to the space of
density matrices is constructed. Further, the projection map thus obtained
is shown to equip the space of regular density matrices with the SLD-Fisher metric,
so that the projection proceeded in \cite{UHI} is a new geometric realization
of the quantum state space (QSS).
\par
A geometric study analogous to \cite{MW} about the projection of the Grover-type
search on the QSS was however deferred to \cite{U} because of another new
interest in the gradient system on the QSS associated with the negative
von-Neumann entropy. It is shown in \cite{UHI} that the gradient system of
interest is understood to be a very natural extension of the gradient system
associated with the negative Shannon entropy on a classical statistical manifold
which is studied by Nakamura \cite{N2}.
This result encourages the author to seek other noted dynamical systems which are
extendable on the QSS. Among the systems displayed in the series of papers
\cite{N1,N2,N3} by Nakamura on integrable systems, the author and Yuya succeed
to extend in \cite{UY1} the averaged Hebbian learning equation (AHLE)
which describes the Hebbian synaptic-neuron learning model proposed by Oja
\cite{N1,Oja}. The dynamical system thus extended on the QSS from the AHLE is
the EAHLE dealt with in this paper. A continuous-time limit of Karmarkar's
projective scaling algorithm of non-constraint \cite{N3,Karmarkar} is also shown
by the author and Yuya to be extendable on the QSS \cite{UY2}.
\par
After the papers \cite{UY1,UY2} by the author and Yuya, the deferred task for a
geometric study of the projection of the Grover-type search sequence for an ordered
tuple of multi-qubits is made successfully by the author \cite{U}: The projection of
the Grover-type search sequence on the QSS is shown to be on an m-geodesic,
an autoparallel curve with respect to the mixture-type parallel transport \cite{Hayashi},
on the QSS. This result strongly encourages the author to seek other dynamical systems
whose trajectories realize geodesics on the QSS.
\par
The aim of this paper is to study the extended averaged Hebbian learning equation
(EAHLE) constructed in \cite{UY1} from a quantum-information-geometric point of
view: All the trajectories of the EAHLE are shown to be the e-geodesics, the
autoparallel curves with respect to the exponential-type parallel transport
\cite{Hayashi}, on the QSS which are known to play an important role not only in
quantum information geometry but also in quantum estimation. In what follows,
the organization of this paper is outlined.
\par
At the beginning of the outline, it should be remarked that the pair of sections,
Section~2 and Section~3, among five sections in this paper are mostly for
reviews of the QSS, the EAHLE and the e-geodesics, which are done according to
the author's previous papers \cite{UY1,UHI,U} and the literature \cite{Hayashi}
by Hayashi on quantum information. Although the review part seems to occupy
rather large part of this paper, it is indispensable because
a similarity between looks of the EAHLE and of the tangent vector along the 
e-geodesics play a key role to reach to the main theorem of this paper. 
\par
In Section~2, the QSS is introduced together with the symmetric logarithmic
derivative (SLD) and the SLD-Fisher metric on that.
The SLD works not only in the definition of the SLD-Fisher metric
endowed with the QSS but also in that of the exponential-type parallel
transport dealt with in Sec.~3. Section~3 is devoted to reviewing
the extended averaged Hebbian learning equation (EAHLE) and the
e-geodesics on the QSS.  The EAHLE is reviewed in subsection~3.1 together with
the way how the EAHLE comes from its origin, the AHLE. In subsection~3.2,
the e-geodesics on the QSS are defined to be the autoparallel curves with respect
to the exponential-type parallel transport. In Section~4, the main theorem of this
paper is proved, which shows that all the trajectories of the EAHLE are the
e-geodesics on the QSS. An explicit representation of solution of the AHLE  is derived
as an outcome of the main theorem. Section~5 is for conclusion.
%%
%%
%%=======================================
%%%%%%%%%%%%%%%%%%%%%%%%%%%%%%%
%
\section{The QSS}
In this section, we set up the quantum state space (QSS) as the space
of regular density matrices endowed with the SLD-Fisher metric,
following \cite{UY1,UHI,UY2}. The literature \cite{Hayashi} by Hayashi is worth
cited to have a general framework on quantum information geometry including
the QSS, in which the QSS is referred to as \lq the space of quantum
states'.
\par
Let $Q_n$ be the set of $n \times n$ regular density matrices, namely, the set
of $n \times n$ positive-definite Hermitean matrices with unit trace. 
On denoting by $M(n)$ the set of all the $n \times n$ complex matrices,
$Q_n$ is defined to be the set
\begin{align}
\label{def-Q}
&
Q_n =
\left\{
\rho \in M(n) \, \left\vert \, \rho:\mbox{positive definite}, \,
\rho^{\dagger}=\rho , \, \mbox{Tr} \,  \rho =1 \, 
\right. \right\}
,
\end{align}
where ${}^{\dagger}$ stands for the Hermitean conjugate operation and
$\mbox{Tr}$ for the trace of matrices. The tangent space, denoted by
$T_{\rho}Q_n$, of $Q_n$ at $\rho \in Q_n$ then takes the form
\begin{align}
\label{def-tangent}
&
T_{\rho}Q_n =
\left\{
X \in M(n) \, \left\vert \, X^{\dagger}=X , \, \mbox{Tr} \, X =0 \, 
\right. \right\}
,
\end{align}
which is equipped with the ${\bf R}$-vector space structure.
\par
As a natural quantum-information-geometric structure of $Q_n$,
the SLD-Fisher metric is endowed with $Q_n$ in what follows.
For the endowment, we need to introduce the symmetric logarithmic
derivative (SLD) on tangent vectors.
The SLD, denoted by $L_{\rho}(X)$, on $X \in T_{\rho}Q_n$ is
the $n \times n$ matrix determined uniquely by the equation
\begin{align}
\label{def-SLD}
&
X
=
\frac{1}{2} \left\{ \, \rho L_{\rho}(X) + L_{\rho}(X) \, \rho \right\}
\quad
(X \in T_{\rho}Q_n ) .
\end{align}
It follows from (\ref{def-SLD}) that the SLD satisfies
\begin{align}
\label{SLD-Hermitean}
&
\left( L_{\rho}(X) \right)^{\dagger}
=
L_{\rho}(X)
\quad
(X \in T_{\rho}Q_n ).
\end{align}
\par
The matrix-element display of the SLD given in \cite{UY1,UHI,UY2} is
of great help also in this paper, which is utilized to prove our main theorem
in Sec.4. Let $\rho \in Q_n$ be written in the form
\begin{align}
\label{diag-rho}
&
\rho 
=
h \, \mathrm{diag} \, 
\left(  \theta_1 , \theta_2 \cdots , \theta_n  \right)
\, h^{\dagger} 
\quad
(h \in U(n)),
\end{align}
where $\mathrm{diag} \, ( \theta_1 , \theta_2 \cdots , \theta_n)$
denotes the $n \times n$ diagonal matrix whose $j$-th diagonal
entry is $\theta_j$ ($j=1,2,\cdots,n$) and  $U(n)$ stands for the group
of $n \times n$ unitary matrices.
The symbol \lq $\mathrm{diag}$' indicates the diagonal matrices
henceforce. On expressing $X \in T_{\rho}Q_n$ as
\begin{align}
\label{X}
&
X=h \wt{X} h^{\dagger}
\end{align}
with $h \in U(n)$ of (\ref{diag-rho}),  
the $(j,k)$-entry of $h^{\dagger} L_{\rho}(X) h$ ($X \in T_{\rho}Q_n$)
is calculated to be
\begin{align}
\label{alt-SLD}
&
\left( h^{\dagger} L_{\rho}(X) h  \right)_{jk}
=
\left( \frac{2}{\theta_j + \theta_k} \right)  \wt{X}_{jk}
\quad
(j,k=1,2, \cdots ,n).
\end{align}
Equations (\ref{def-SLD})-(\ref{alt-SLD}) are put together to
show the following lemma on the SLD.
\begin{lemma}
\label{lemma_SLD}
The symmetric logarithmic derivative (SLD), $L_{\rho}$, is a one-to-one and
onto ${\bf R}$-linear map from $T_{\rho}Q_n$ to
\begin{align}
\label{SLD-Q_n}
&
L_{\rho}(T_{\rho} Q_n)
=
\left\{
\Xi \in M(n) \, \left\vert \, \, \Xi^{\dagger}=\Xi , \, 
\Tr (\rho \Xi + \Xi \rho) =0 \right.
\right\} 
\quad
(\rho \in Q_n) . 
\end{align}
Under (\ref{diag-rho}) and
\begin{align}
\label{Xi}
&
\Xi = h \wt{\Xi} h^{\dagger} \in L_{\rho}(T_{\rho}Q_n) ,
\end{align}
the inverse, denoted by $L_{\rho}^{-1}$, of the SLD
is given to satisfy
\begin{align}
\label{inv-SLD}
&
\left(  h^{\dagger} L_{\rho}^{-1}(\Xi) h  \right)_{jk}
=
\left( \frac{\theta_j+\theta_k}{2} \right) \wt{\Xi}_{jk}
\quad (j,k=1,2,\cdots ,n).
\end{align}
\end{lemma}
\par
In terms of the SLD, the SLD-Fisher metric $\langle \cdot , \cdot \rangle$
is defined by
\begin{align}
\label{def-Fisher}
&
\langle X , \, X^{\prime} \rangle_{\rho}
=
\Tr \left( X^{\dagger} L_{\rho}(X^{\prime}) \right)
\quad
(X, X^{\prime} \in T_{\rho}Q_n)
\end{align}
(see Hayashi \cite{Hayashi}). On using (\ref{def-SLD}) and
(\ref{SLD-Hermitean}), the SLD-Fisher metric is brought into
the form
\begin{align}
\label{SLD-Fisher}
&
\langle X , \, X^{\prime} \rangle_{\rho}
=
\frac{1}{2} \,
\Tr
 \left(
        \rho 
         \left( L_{\rho}(X) L_{\rho}(X^{\prime})
                 +  L_{\rho}(X^{\prime}) L_{\rho}(X)  \right)
\right)
\quad
(X, X^{\prime} \in T_{\rho}Q_n) 
\end{align}
\cite{UY1,UHI,UY2,U}. Furthermore, with the matrix-element displays,
(\ref{diag-rho})-(\ref{alt-SLD}), and
\begin{align}
\label{X'}
&
X^{\prime}=h \wt{X}^{\prime}  h^{\dagger},
\end{align}
the SLD-Fisher metric is expressed to be
\begin{align}
\label{alt-Fisher}
&
\langle X , X^{\prime} \rangle_{\rho}
=
\sum_{j,k=1}^{n}
\left( \frac{2}{\theta_j+\theta_k} \right) \,
\overline{\wt{X}}_{jk} \wt{X}_{jk}^{\prime}
 \end{align}
\cite{UY1,UHI,UY2,U}. Equation (\ref{alt-Fisher}) works effectively
to derive a useful formula to the gradient equation on the QSS
\cite{UY1,UHI,UY2}.
The Riemannian manifold $Q_n$ endowed with the SLD-Fisher metric
$\langle \cdot , \cdot \rangle$ is what we are referring to as the quantum
state space.
%%
%%
%%============================================
%%%%%%%%%%%%%%%%%%%%%%%
\section{The EAHLE and the e-geodesics}
%}
In this section, the extended averaged Hebbian learning equation (EAHLE)
and the e-geodesics are introduced according to \cite{UY2} for the EAHLE
and to \cite{Hayashi} for the e-geodesics.
\subsection{The EAHLE}
The extended averaged Hebbian learning equation (EAHLE) is
organized by the present author and Yuya \cite{UY2} who are
inspired by Nakamura's paper \cite{N1} on the averaged Hebbian learning
equation (AHLE). According to \cite{UY2}, the EAHLE
is the first order differential equation
\begin{align}
\label{EAHLE}
&
\frac{d\rho}{dt}
=
\rho \, C + C \rho  -2 \, \Tr (C\rho) \, \rho
\end{align}
on the QSS. The $C$ on the rhs of (\ref{EAHLE}) is the real diagonal matrix
\begin{align}
\label{C}
&
C= \mathrm{diag} \, \left(  c_1, c_2, \cdots , c_n  \right)
\end{align}
of degree $n$, whose diagonal entries stand for the eigenvalues
of the autocorrelation matrix of the stationary stochastic process
governing the Hebbian learning process \cite{N1,Oja,Hebb}.
\par
The reason for referring to Eq.\ (\ref{EAHLE}) as the \lq extended'
averaged Hebbian learning equation is given in what follows.
Let $w=(w_j)_{j=1,2,\cdots ,n}$ be the variables on $t$
expressing the coupling strengths of neurons obtained through an
appropriate change of the independent variable to have $t$.
The first-order differential equation
\begin{align}
\label{AHLE}
&
\frac{dw} {dt}
=
Cw - \left(  w^T C w  \right) \, w
\quad
\end{align}
on the $n-1$ dimensional sphere
\begin{align}
\label{sphere}
&
S^{n-1}
=
\left\{  w = (w_1,w_2, \cdots ,w_n)^T \in {\bf R}^n 
\, \left\vert \,
w^T w =1  \right. 
\right\} 
\end{align}
with unit radius describes Oja's rule \cite{Oja} on the Hebbian learning process
of synaptic neurons \cite{Hebb}, where $C$ is the diagonal matrix given by
(\ref{C}). In (\ref{AHLE}), the $C$ is understood again to be the diagonalization
of the autocorrelation matrix of the governing stationary stochastic process
of neurons. The differential equation (\ref{AHLE}) is what we are referring
to as the averaged Hebbian learning equation (AHLE).
\par
If we restrict Eq.\ (\ref{AHLE}) on each of the open subsets,
\begin{align}
\nonumber
&
S^{n-1}_{\sigma}
=
\left\{
w \in S^{n-1} \, \left\vert \, \sigma_j w_j >0, \, j=1,2,\cdots ,n 
\right. \right\}
\\
\label{S_sigma}
&
\quad
\left(
\sigma=(\sigma_j), \, \sigma_j =\pm 1 , \, j=1,2,\cdots ,n
\right)
,
\end{align}
of $S^{n-1}$, it is brought into the first-order differential equation
\begin{align}
\label{AHLE-D}
&
\frac{ d\theta_j}  {dt}
=
2c_j \theta_j - 2 \left( \sum_{k=1}^{n} c_j \theta_j \right) \, \theta_j
\quad
(j=1,2,\cdots,n) ,
\end{align}
on the submanifold
\begin{align}
\label{D_n}
&
D_n 
=
\left\{
\Theta \in Q_n \, \left\vert \, 
\Theta = \mathrm{diag} \, (\theta_1,\theta_2, \cdots , \theta_n) \,
\right.
\right\} 
\end{align}
of $Q_n$ through the map
\begin{align}
\nonumber
&
p_{n,\sigma} (w) = 
\mathrm{diag} \left(  w_1^2,  w_2^2,\cdots , w_n^2  \right) 
\\
\label{p_n_sigma}
&
\quad
\left(
w \in S^{n-1}_{\sigma}, \,
\sigma=(\sigma_j), \, \sigma_j =\pm 1 , \, j=1,2,\cdots ,n
\right)
\end{align}
from $S^{n-1}_{\sigma}$ to $D_n$ \cite{UY1}. 
The $\theta_j$'s in (\ref{AHLE-D}) and (\ref{D_n}) are subject to the
constraints 
\begin{align}
\label{theta-constraints}
&
\theta_j > 0 \,\, (j=1,2,\cdots n)
\quad \mbox{and} \quad
\sum_{j=1}^{n} \theta_j =1 .
\end{align}
We note here that Eq.\ (\ref{AHLE-D}) is the same form
as the Toda lattice written in Moser's form \cite{N1,Moser,MZ}.
\par
Since every mapping $p_{n,\sigma}$ defined by  (\ref{p_n_sigma}) is
a diffeomorphism with the inverse
\begin{align}
\nonumber
&
p_{n,\sigma}^{-1} (\Theta)
=
\left(
\sigma_1 \sqrt{\theta_1}, \sigma_2 \sqrt{\theta_2},
\cdots , \sigma_n \sqrt{\theta_n}
\right)
\\
\label{p_n_sigma_inv}
&
\quad
\left(
\Theta \in D_n, \,
\sigma=(\sigma_j), \, \sigma_j =\pm 1 , \, j=1,2,\cdots ,n
\right)
,
\end{align}
we can understand that the differential equation (\ref{AHLE-D})
on $D_n$ is a \lq copy' of the AHLE restricted on each $S^n_{\sigma}$
and vice versa \cite{UY1}.
\par
We are at the final stage to account for the naming of the EAHLE.
To complete the account, we show that the differential equation
(\ref{AHLE-D}) on $D_n$ is the restriction of the EAHLE, (\ref{EAHLE}),
on $D_n$. In fact, the substitution of $\Theta \in D_n$ for $\rho$ in (\ref{EAHLE})
yields the \lq copy', (\ref{AHLE-D}), of the AHLE.
In a summary, we have the following lemma \cite{UY1,N1}.
\begin{lemma}
\label{AHLE-EAHLE}
The restriction of the averaged Hebbian learning equation (AHLE) on each
$S^{n-1}_{\sigma}$ is equivalent,
up to the diffeomorphisms given by (\ref{p_n_sigma}), 
to Eq.\ (\ref{AHLE-D}) which describes not only the restriction of the extended
averaged Hebbian learning equation (EAHLE) on $D_n$ but also the Toda lattice
in Moser's form \cite{N1,Moser,MZ}.
\end{lemma}
%%%%
%%%%
\subsection{The e-geodesics on the QSS}
To those who are not familiar with differential geometry,
a geodesic connecting a given pair of points would be thought of
as the shortest-distance path between the given points. For example,
in the Euclidean plane, a typical model space
for school-geometry, we are taught that the straight-line segment
connecting a given pair of points is the geodesic between them.
In differential geometry, however, the notion of length or distance,
is unnecessary in defining geodesics:
What is needed in the definition of geodesics is the idea of
parallel transports, namely, the idea
for comparing tangent vectors at a certain point
with those at another point. Once a parallel transport is fixed,
the geodesics are defined to be the autoparallel curves with respect
to that parallel transport.
For an intuitive description of geodesics and parallel transports, 
the literature \cite{Nakahara} by Nakahara is worth cited.
\par
Let us start with the definition
of the exponential-type (e-) parallel transport.
According to Hayashi \cite{Hayashi},
the e-parallel transport from $T_{\rho_1}Q_n$ to
$T_{\rho_2}Q_n$ is the ${\bf R}$-linear map
$\tau_{\rho_1,\rho_2}: T_{\rho_1}Q_n \rightarrow T_{\rho_2}Q_n$
subject to
\begin{align}
\label{parallel}
&
L_{\rho_2}(\tau_{\rho_1,\rho_2}(X))
=
L_{\rho_1}(X)- \Tr \left( \rho_2 \, L_{\rho_1}(X) \right)
\quad
(X \in T_{\rho_1}Q_n) ,
\end{align}
where $L_{\rho_1}$ and $L_{\rho_2}$ denote the SLD defined
by (\ref{def-SLD}) with $\rho =\rho_1$ and $\rho=\rho_2$, respectively.
Combining the defining equation (\ref{def-SLD}) for the SLD
with Eq.\ (\ref{parallel}), we can obtain a more direct form,
\begin{align}
\label{alt-parallel}
&
\tau_{\rho_1,\rho_2}(X)
=
\frac{1}{2}\left\{ \rho_2  L_{\rho_1}(X) + L_{\rho_1}(X) \rho_2 \right\}
-
\Tr \left( \rho_2 \, L_{\rho_1}(X) \right)  \rho_2 
\quad
(X \in T_{\rho_1}Q_n) ,
\end{align}
of the e-parallel transport. The e-parallel transport satisfies, of course,
the postulate of parallel transports (see Guggenheimer \cite{Gugg},
for example), but we do not get it into detail here. 
\begin{definition}
\label{def-e-parallel}
Tangent vectors $X_1 \in T_{\rho}Q_n$ and $X_2 \in T_{\rho_2} Q_n$
are e-parallel if they are parallel with respect to the e-parallel transport;
namely, $X_1$ and $X_2$ are e-parallel if they satisfy
\begin{align}
\label{def-para}
&
X_2 = \tau_{\rho_1,\rho_2}(X_1) ,
\end{align}
where $\tau_{\rho_1, \rho_2}$ is the e-parallel transport given by
(\ref{alt-parallel}).
\end{definition}
\par
Once we fix a parallel transport, we can consider the geodesics
as the autoparallel curves to that parallel transport \cite{Nakahara,Gugg}.
In view of Definition \ref{def-e-parallel},
we can define the e-geodesics as follows \cite{Hayashi}.
\begin{definition}
\label{def-e-geod}
A smooth curve $\rho (t)$ ($0 \leq t \leq {}^{\exists}T$) 
on the QSS is an e-geodesic if it satisfies
\begin{align}
\nonumber
\frac{d \rho}{dt}(t) 
&= \tau_{\rho (0) , \, \rho (t)} \left( \frac{d \rho}{dt}(0) \right)
\\
\nonumber
&=
\frac{1}{2} 
\left\{ \rho (t) L_{\rho (0) } \left( \frac{d\rho}{dt}(0) \right) 
+ L_{\rho (0) } \left( \frac{d\rho}{dt}(0) \right) \rho (t) \right\}
\\
\label{autoparallel}
& \quad \qquad
-\Tr \left(
L_{\rho (0)} \left( \frac{d\rho}{dt}(0)) \right) \rho (t) 
\right)  \rho (t)
\qquad
(0 \leq t \leq T) ,
\end{align}
where $\tau_{\rho (0) , \, \rho (t)}$ is the e-parallel transport
from $T_{\rho(0)}Q_n$ to $T_{\rho (t)}Q_n$ given by (\ref{alt-parallel})
with $\rho (0)$ and $\rho (t)$ in place of $\rho_1$ and $\rho_2$, respectively.
\end{definition}
Equation (\ref{autoparallel}) for the autoparallelism with respect to
the e-parallel transport (\ref{alt-parallel}) must not be understood to be
a first-order differential equation on the QSS because of the appearance
of the initial tangent vector $(d\rho /dt)(0)$ in the rhs of (\ref{autoparallel})
that never takes place in first-order differential equations.
Hence the expression (\ref{autoparallel}) of the e-geodesics does not
contradict the second-order-differential-equation form taught in theory
of geodesics. According to Hayashi \cite{Hayashi}, the e-geodesic admits
the explicit representation below.
\begin{lemma}
\label{e-geod-exp}
The e-geodesic $\rhoe \left( t; \rhozero , \Xzero \right)$ with the initial conditions,
\begin{align}
\label{init-pos}
&
\rhoe \left( 0;\rhozero, \Xzero \right) = \rhozero \in Q_n
\\
\intertext{and}
&
\label{init-tan}
\frac{d\rhoe}{dt} \left( 0; \rhozero,\Xzero \right) 
= \Xzero \in T_{\rhozero} Q_n ,
\end{align}
takes the form
\begin{align}
\nonumber
\rhoe \left( t; \rhozero, \Xzero \right)
=
&
\left\{
\Tr \left( e^{\frac{t}{2}L_{\rhozero}(\Xzero)} \rhozero \, e^{\frac{t}{2}L_{\rhozero}(\Xzero)} \right)
\right\}^{-1}
\\
\label{e-geod}
&
\quad
\times 
e^{\frac{t}{2}L_{\rhozero}(\Xzero)} \rhozero \, e^{\frac{t}{2}L_{\rhozero}(\Xzero)}
\end{align}
for $0 \leq t < \infty$, where $L_{\rhozero} \left( \Xzero \right)$
is the SLD, defined by (\ref{def-SLD}), on $\Xzero \in T_{\rhozero}Q_n$.
\end{lemma}
A direct differentiation of (\ref{e-geod})  by $t$ clearly shows that
$\rhoe (t;\rhozero , \Xzero )$ satisfies Eq.\ (\ref{autoparallel})
of the autoparallelism. Namely, we have
\begin{align}
\nonumber
&
\frac{d \rhoe}{dt} \left( t; \rhozero,\Xzero \right)
\\
\nonumber
=
&
\frac{1}{2} 
\left\{
\rhoe \left( t;\rhozero,\Xzero \right)
L_{\rhozero} \left( \Xzero \right) 
+
L_{\rhozero} \left( \Xzero \right)
\rhoe \left( t;\rhozero, \Xzero \right)
\right\}
\\
\label{e-geod-diff}
&
\quad
-\Tr \left(
L_{\rhozero}\left( \Xzero \right) \rhoe \left (t;\rhozero,\Xzero \right) \right)
\, \rhoe \left( t;\rhozero,\Xzero \right) .
\end{align}
We see that Eq.\ (\ref{e-geod-diff})
for any fixed e-geodesic looks quite similar to the EAHLE (\ref{EAHLE}) 
and therefore we may expect that every trajectory of
the EAHLE can be realized as an e-geodesic.
%%
%%
%============================================================
%%
%%
%%%%%%%%%%%%%%%%%%%%%%%%%%%%%%%%%%
\section{The EAHLE trajectories as the e-geodesics}
Now that we have found the similarity of the EAHLE (\ref{EAHLE})
and Eq.\ (\ref{e-geod-diff}) for the e-geodesics, we show that
the trajectories of the EAHLE are the e-geodesics below.
Further, as an outcome of Theorem~\ref{main}, an explicit representation
of solution of the AHLE is derived from the representation (\ref{e-geod})
of the e-geodesics.
\par
On comparing very naively Eq.\ (\ref{EAHLE}) with
Eq.\ (\ref{e-geod-diff}), one might come to choose
$L_{\rhozero }( \Xzero )$ in (\ref{e-geod-diff}) to be equal to
$C$ in (\ref{EAHLE}).
However, this choice fails because the diagonal matrix $C$ never
belongs to $L_{\rhozero} (T_{\rhozero}Q_n )$
(see (\ref{SLD-Q_n}) with $\rhozero$ in place of $\rho$).
By choosing $L_{\rhozero}(\Xzero)$  in (\ref{e-geod-diff}) suitably,
we have the following theorem to characterize all the trajectories
of the EAHLE as the e-geodesics.
\par
\begin{theorem}[Main Theorem]
\label{main}
For any fixed $\rhozero \in Q_n$, 
let $\rhoh (t; \rhozero)$ denote the trajectory of the EAHLE
subject to the initial condition
%%
%\smallskip
%%
\begin{eqnarray}
\label{init-EAHLE}
\rho^h (0;\rhozero)=\rhozero , 
\end{eqnarray}
%%
%\smallskip
%%
and let $\rhoe (t; \rhozero, \XC)$ denote
the e-geodesic $\rhoe (t; \rhozero, \XC)$ subject to the initial conditions
\begin{align}
\label{init-e-geod-pos}
&
\rhoe (0; \rhozero, \XC) = \rhozero
\\
\intertext{and}
\label{init-e-geod-tan}
&
\frac{d\rhoe}{dt} (0; \rhozero, \XC) = \XC
\\
\intertext{with}
\label{XC}
&
\XC
=
\rhozero C + C \rhozero - 2\mathrm{Tr}\,(C \rhozero) \rhozero .
\end{align}
Then, for $t \geq 0$,
the trajectory $\rhoh (t; \rhozero)$ of the EAHLE coincides with
the e-geodesic $\rhoe (t; \rhozero, \XC)$.
\end{theorem} 
\par\noindent
{\bf Proof}
\quad
%\par\noindent
%%
From (\ref{init-EAHLE}) and (\ref{init-e-geod-pos}),
we easily confirm that $\rhoh (t; \rhozero )$ and
$\rhoe (t;\rhozero, \XC)$ share $\rhozero$ as the initial point.
As a necessary condition for
$
\rhoe (t; \rhozero , \XC) = \rhoh (t; \rhozero )
$
,
we pose the coincidence
%$
\begin{align}
\label{init-tan-coincidence}
&
\frac{d\rhoe}{dt}(0; \rhozero , \XC)
=
\frac{d\rhoh}{dt}(0; \rhozero )
\end{align}
of the initial tangent vectors, 
which is equivalent to the pair of equations, (\ref{init-e-geod-tan}) and (\ref{XC}).
Hence, what we have to show here is that the pair,
 (\ref{init-e-geod-tan}) and (\ref{XC}), is also a sufficient condition for
$
\rhoe (t; \rhozero , \XC) = \rhoh (t; \rhozero )
$
.
Through the proof, we often apply the abbreviation $\rho^e$ to
$\rho^e(t:\rhozero , \XC )$.
%%
%\par
%%
The core part of the proof is given in what follows in a straightforward
calculation form with the notation
\begin{align}
\label{init-diag}
&
h^{\dagger} \rhozero h
=
\Thetaz
=
\mathrm{diag} \, (\thetaz_1  , \thetaz_2 , \cdots , \thetaz_n)
\quad
(h \in U(n)) 
\\
\intertext{and}
\label{notation}
&
\wt{X}\supC = h^{\dagger} \XC h , \quad
\wt{L}\supz = h^{\dagger} L_{\rhozero}(\XC) h , \quad
\wt{C} = h^{\dagger} C h , \quad
R= h^{\dagger} \rho^e h .
\end{align}
\noindent
We note here that $\wt{C}$ and $R$ are non-diagonal in general
and that $R$ is of unit trace. 
With the notation given above,
the matrix-element display of $\wt{X}\supC$ takes the form
\begin{align}
\nonumber
\wt{X} \supC _{jk}
&=
\left(
h^{\dagger} \{ \rhozero C + C \rhozero - 2\, \Tr (C\rhozero) \rhozero \}h
\right)_{jk}
\\
\nonumber
&=
\left(
 \Thetaz \wt{C} + \wt{C} \Thetaz 
 -
 2\, \Tr (\wt{C} \Thetaz ) \Thetaz \right)_{jk}
\\
\label{X^C_jk}
&=
\thetaz _j \wt{C}_{jk} + \wt{C}_{jk} \thetaz _k
- 2 \, \Tr (\wt{C} \Thetaz ) \delta_{jk} \thetaz _k 
\quad
(j,k=1,2,\cdots,n) ,
\end{align}
where the symbol $\delta_{jk}$ indicates Kronecker's delta
($j,k=1,2,\cdots ,n$).
Equation (\ref{X^C_jk}) is combined with (\ref{alt-SLD}) to yield
the matrix-element display
\begin{align}
\label{tilde-L_jk}
\wt{L}\supz _{jk}
=
2 \wt{C}_{jk} 
-
\left(
\frac{4}{\thetaz _j + \thetaz _k}
\right)
\,
\Tr \left( \wt{C} \Thetaz \right) \delta_{jk} \thetaz _k 
\quad
(j,k=1,2,\cdots,n) 
\end{align}
for $\wt{L}\supz$.
\par
We are now in a position to calculate the rhs of (\ref{e-geod-diff})
with $\Xzero = \XC$ and (\ref{XC}).
Putting Eqs.\ (\ref{init-diag}), (\ref{notation}) and (\ref{tilde-L_jk})
together with the abbreviation $\rho^e$ for $\rhoe (t; \rhozero, \XC)$,
we have
\begin{align}
\nonumber
&
\left(
h^{\dagger} 
\left\{
\frac{1}{2}
\left(
\rho^e L_{\rhozero} \left( \XC \right) + L_{\rhozero}\left( \XC \right) \rho^e 
\right)
- \Tr \left( L_{\rhozero}(\XC) \rho^e \right) \rho^e 
\right\} h 
\right)_{jk}
\\
\nonumber
=
&
\left(
\frac{1}{2}
\left(
R \wt{L}\supz + \wt{L}\supz R
\right) 
- \Tr \left( \wt{L}\supz R \right) R 
\right)_{jk}
\\
\nonumber
=
&
\frac{1}{2} \sum_{m=1}^n 
\left\{
\left( R_{jm} \wt{L}\supz _{mk} + \wt{L}\supz _{jm} R_{mk} \right) 
- 2 \Tr (\wt{L}\supz R) R_{jk} \right\}
\\
\nonumber
=
&
\sum_{m=1}^n
R_{jm} 
\left\{
\wt{C}_{mk} - \left( \frac{2 }{\thetaz_m + \thetaz_k} \right)
   \Tr \left( \wt{C}\Thetaz \right) \delta_{mk} \thetaz_k
\right\}
\\
\nonumber
&
\qquad
-
\sum_{m=1}^n
\left\{
\wt{C}_{jm} - \left( \frac{2}{\thetaz_j + \thetaz_m} \right) 
\Tr \left( \wt{C}\Thetaz \right) \delta_{jm} \thetaz_m \right\} 
R_{mk}
\\
\nonumber
&
\qquad \qquad
-
2
\left[
\sum_{m,l=1}^{n}  
\left\{
\wt{C}_{lm} 
- \left( \frac{2}{\thetaz_l + \thetaz_m} \right) \Tr \left( \wt{C}\Thetaz \right)
\delta_{lm} \thetaz_m \right\} R_{ml} 
\right]
R_{jk}
\\
\noalign{\bigskip}
\nonumber
=
&
\sum_{m=1}^{n} 
\left( R_{jm}\wt{C}_{mk} + \wt{C}_{jm}R_{mk} \right)
-2 \Tr \left( \wt{C}\Thetaz \right) R_{jk} 
\\
\nonumber
&
\qquad
-2 \left( \sum_{m,l=1}^n \wt{C}_{lm}R_{ml} \right) R_{jk}
+ 2 \Tr \left( \wt{C}\Thetaz \right) \left( \sum_{m=1}^n R_{mm} \right) R_{jk}
\\
\nonumber
=
&
\left( R \wt{C} + \wt{C}R \right)_{jk} - 2 \Tr \left( \wt{C}R \right) R_{jk}
\\
\label{rhs-e-geod-diff}
=
&
\left(
h^{\dagger}
\left\{ \rho^e C + \rho^e C - 2\Tr (C\rho^e) \rho^e \right\}
h
\right)_{jk} .
\end{align}
Equation (\ref{rhs-e-geod-diff}) is put together with
Eq.\ (\ref{e-geod-diff}) to show that $\rho^e(t:\rhozero , \XC )$
satisfies the equation
\begin{align}
\nonumber
\frac{d\rhoe}{dt}(t:\rhozero , \XC )
=
&
\rhoe \left( t:\rhozero , \XC \right) C
+
C \rhoe \left( t:\rhozero , \XC \right)
\\
\label{e-geod-EAHLE}
&
\quad
-2\Tr (C \rho^e(t:\rhozero , \XC )) \rho^e(t:\rhozero , \XC ) .
\end{align}
which turns out to be the same as Eq.\ (\ref{EAHLE}) with
$\rhoe  (t;\rhozero,\XC)$ in place of $\rho$. Put in another way,
the e-geodesic $\rhoe (t;\rhozero, \XC)$ satisfies the
EAHLE (\ref{EAHLE}) with the initial condition $\rho (0)=\rhozero$
and accordingly $\rhoe (t;\rhozero, \XC)$ coincides with
$\rhoh (t;\rhozero)$.
This completes the proof.
\par\medskip
Combining Lemma~\ref{AHLE-EAHLE}, Lemma~\ref{e-geod-exp} and
Theorem~\ref{main} together, we can give an explicit representation of
solution of the averaged Hebbian learning equation (AHLE) from
the representation, (\ref{e-geod}), of the e-geodesics.
Although it is pointed out in Nakamura \cite{N1} that the representation
of solution of the AHLE is available from that of the Toda lattice in Moser's
form \cite{Moser,MZ} (see also Lemma~\ref{AHLE-EAHLE}),
we are to present the same one as in \cite{N1} because our derivation
process below is new.
\par
On recalling Lemma~\ref{AHLE-EAHLE}, the representation of solution
of the AHLE can be given by calculating explicitly the Eq.\ (\ref{e-geod})
with the initial conditions,
\begin{align}
\label{init-diag-pos}
&
\rhoe (0) = \rhozero = \Thetaz 
= \diag
   \left(\thetaz_1    , \thetaz_2 , \cdots , \thetaz_n \right) \in D_n \subset Q_n
\\
\intertext{and}
\label{init-diag-velocity}
&
\frac{d\rhoe}{dt}(0)
=
\Xzero
=
\XiC
=
\diag \left( \xiz_1    , \xiz_2 , \cdots , \xiz_n \right) \in T_{\Thetaz}D_n
\subset T_{\Thetaz}Q_n
\\
\intertext{with}
\label{xiz}
&
\xiz_j
=
2c_j \thetaz_j -2 \mathrm{Tr} \, \left( C\Thetaz \right) \, \thetaz_j 
\qquad
(j=1,2,\cdots,n),
\end{align}
where $C$ is the diagonal matrix governing both the EAHLE (\ref{EAHLE})
and the AHLE (\ref{AHLE}).
Under (\ref{init-diag-pos})-(\ref{xiz}), the SLD of $\XiC$ and its exponential
are calculated to be
\begin{align}
\label{SLD-Xc-diag}
&
L_{\Thetaz}\left( \XiC \right)
=
\diag 
\left(
\frac{\xiz_1}{\thetaz_1}, \frac{\xiz_2}{\thetaz_2}, \cdots ,
\frac{\xiz_n}{\thetaz_n}
\right)
=
2C- 2 \mathrm{Tr} \, \left( C\Thetaz \right) I 
\\
\intertext{and}
&
\label{exp-Thetaz}
e^{\frac{t}{2}L_{\Thetaz} \left( \XiC \right)}
=
e^{-t \, \mathrm{Tr} \, \left( C\Thetaz \right)} \,
\diag
\left(
e^{tc_1}, e^{tc_2}, \cdots , e^{tc_n}
\right)  .
\end{align}
Hence it follows from Lemma~\ref{e-geod-exp} and Theorem~\ref{main}
that the trajectory on $D_n$, denoted by $\Thetah (t)$,
of the EAHLE with the initial condition (\ref{init-diag-pos}) takes the form
\begin{align}
\label{EAHLE-diag-1}
&
\Thetah (t) 
=
\diag \left( \thetah_1(t), \thetah_2(t), \cdots , \thetah_n(t) \right)
\\
\intertext{with}
&
\label{EAHLE-diag-2}
\thetah_j (t)
=
\left( \sum_{k=1}^n  e^{2tc_k}\thetaz_k \right)^{-1}
e^{2tc_1}\thetaz_j
\quad
(j=1,2,\cdots,n) .
\end{align}
We note here that Eqs.\ (\ref{EAHLE-diag-1})
and (\ref{EAHLE-diag-2}) reproduce the solution of the Toda lattice
in Moser's form \cite{Moser,MZ} in view of Lemma~\ref{AHLE-EAHLE}.
Applying the map $p_{n,\sigma}^{-1}$ defined by
(\ref{p_n_sigma_inv}) to $\Thetah (t)$, we have
\begin{align}
\nonumber
w_{\sigma}^h(t)
&=
p_{n,\sigma}^{-1}(\Theta^h(t))
\\
\nonumber
&
=
\left( \sum_{k=1}^n  e^{2tc_k} \thetaz_k \right)^{-1/2}
\left(
e^{tc_1}\sigma_1 \sqrt{\thetaz_1},
e^{tc_2}\sigma_2 \sqrt{\thetaz_2},
\cdots ,
e^{tc_n}\sigma_n \sqrt{\thetaz_n}
\right)^T 
\\
\label{inv-Theta}
& \qquad
\left( \sigma=(\sigma_j), \, \sigma_j =\pm 1 , \, j=1,2,\cdots ,n \right)
\end{align}
which realizes the solutions of the AHLE on the open-dense subset
$\cup_{\sigma} S^{n-1}_{\sigma}$ of $S^{n-1}$
(see (\ref{S_sigma}) for $S^{n-1}_{\sigma}$).  
Thus we have the following corollary to Theorem~\ref{main}.
\begin{corollary}
\label{AHLE-sol}
The solution of the averaged Hebbian learning equation (\ref{AHLE}) subject
to the initial condition $w(0)=(\wz_1,\wz_2,\cdots,\wz_n)^T \in S^{n-1}$
is given by
\begin{align}
\label{wh1}
&
\wh(t)=
\left( \sum_{k=1}^n  e^{2tc_k} \left( \wz_k \right)^2 \right)^{-1/2}
\left(
e^{tc_1}\wz_1, e^{tc_2}\wz_j, \cdots , e^{tc_n}\wz_n
\right)^T .
\end{align}
\end{corollary}
%%
%%
%============================================================
%%
%%
%%%%%%%%%%%%%%%%%%%%%%%%%%%%%
\section{Conclusions}
We show in Theorem~\ref{main}
that all the trajectories of the extended averaged
Hebbian learning equation (EAHLE) on the QSS are the e-geodesics
on the QSS. As a direct application of Theorem~\ref{main},
the explicit representation,  (\ref{wh1}), of solution of the averaged
Hebbian learning equation (AHLE), the departure equation
of the EAHLE, is derived from the representation,  (\ref{e-geod}),
of the e-geodesics.
Although the expression (\ref{wh1}) is known already to be
available from the Toda lattice in Moser's form due to the
equivalence between the AHLE and the Moser's form \cite{N1},
our derivation in Sec.4 is worth given because it is made from a novel
point of view, a quantum-information-geometric point of view.
\par
We would like to offer a remark on Theorem~\ref{main} from a
geometric-mechanics point of view: In view of the gradient-system
structure of the EAHLE revealed in Uwano and Yuya \cite{UY1},
Theorem~\ref{main} is understood to provide a gradient system
whose trajectories are the e-geodesics. Further, since the e-geodesics
are known to play an important role in quantum estimation \cite{Hayashi},
the EAHLE is expected to be a new candidate of gradient
systems dealt with in Braunstein \cite{Br} for quantum estimation.
Another remark is offered from an integrable-systems point of view:
The EAHLE would be looked on as an extended Toda lattice in Moser's form:
This view is supported from the coincidence given in Lemma~\ref{AHLE-EAHLE}
between the EAHLE restricted on $D_n$ and the Moser's form.
\par
On closing this paper, the author would like to make the following conjecture
on the e-geodesics which do not satisfy the initial conditions,
(\ref{init-e-geod-tan}) and (\ref{XC}), attached to the initial tangent vector.
\par\medskip\noindent
{\bf Conjecture}\,\,
{\it
All the e-geodesics on the QSS are realized as the EAHLE
trajectories up to the adjoint $SU(n)$ actions on the QSS and
the affine transformations of time.
}
\par\medskip\noindent
The conjecture will be investigated soon together with the dynamics on
the QSS described by the EAHLE.
%%
%%
%############################################################
%%
\par\bigskip\noindent
{\bf Acknowlegement}\,\,
The author thanks Professor Yoshimasa Nakamura
at Kyoto University for his valuable remark on the paper \cite{N1}
and his suggestion to include our derivation of the AHLE solution
(Corollary~\ref{AHLE-sol}) in the present paper.
%%
%%
%%%%%%%%%%%%%%%%%%%%%%%%%%%%%

\end{document}